\documentclass[11pt]{article}

 \usepackage{latexsym}
 \usepackage{amsbsy}
 \usepackage{amssymb}
 \usepackage{amsmath}
 \input{amssym.def}
 \input{amssym}
 \usepackage{graphicx}
 \usepackage{epsfig}
 \usepackage[latin1]{inputenc}
 \usepackage{a4wide}

 \usepackage{epstopdf}

 \newtheorem{thm}{Theorem}[section]
 \newtheorem{defin}[thm]{Definition}
 \newtheorem{lem}[thm]{Lemma}
 \newtheorem{prop}[thm]{Proposition}
 \newtheorem{cor}[thm]{Corollary}
 \newtheorem{rem}[thm]{Remark}
 \newtheorem{ex}[thm]{Example}

 \newcommand{\bthm}{\begin{thm}}
 \newcommand{\ethm}{\end{thm}}
 \newcommand{\bd}{\begin{defin}}
 \newcommand{\ed}{\end{defin}}
 \newcommand{\blem}{\begin{lem}}
 \newcommand{\elem}{\end{lem}}
 \newcommand{\bcor}{\begin{cor}}
 \newcommand{\ecor}{\end{cor}}
 \newcommand{\bprop}{\begin{prop}}
 \newcommand{\eprop}{\end{prop}}
 \newcommand{\brem}{\begin{rem} \rm}
 \newcommand{\erem}{\end{rem}}
 \newcommand{\bex}{\begin{ex} \rm}
 \newcommand{\eex}{\end{ex}}

 \newcommand{\pr}{\noindent{\bf Proof. }}
 \newcommand{\ep}{\nolinebreak{\hspace*{\fill}$\Box$ \vspace*{0.25cm}}}

 \newcommand{\beq}{\begin{equation}}
 \newcommand{\eeq}{\end{equation}}

 \newcommand{\bea}{\begin{eqnarray}}
 \newcommand{\eea}{\end{eqnarray}}

 \newcommand{\beas}{\begin{eqnarray*}}
 \newcommand{\eeas}{\end{eqnarray*}}

 \newcommand{\beqs}{\begin{equation*}}
 \newcommand{\eeqs}{\end{equation*}}

 \newcommand{\bi}{\begin{itemize}}
 \newcommand{\ei}{\end{itemize}}

 \newcommand{\ben}{\begin{enumerate}}
 \newcommand{\een}{\end{enumerate}}

 \newcommand{\ba}{\begin{array}}
 \newcommand{\ea}{\end{array}}

 \newcommand{\R}{\mathbb R}

 \newcommand{\cF}{\ensuremath{{\cal F}}}

 \newcommand{\cL}{\ensuremath{{\cal L}}}

 \newcommand{\eps}{\varepsilon}
 \newcommand{\vphi}{\varphi}

 \newcommand{\rep}{\mathop{\rm Re}\nolimits}
 \newcommand{\imp}{\mathop{\rm Im}\nolimits}

 \newcommand{\tg}{\mathop{\rm tg}\nolimits}
 \newcommand{\ctg}{\mathop{\rm ctg}\nolimits}

 \newcommand{\tgh}{\mathop{\rm tgh}\nolimits}

 \newcommand{\notmid}{\mid\kern-0.5em\not\kern0.5em}

 \newcommand{\pd}{\ensuremath{\partial}}

 \def\dj{d\kern-0.4em\char"16\kern-0.1em}
 \def\Dj{\mbox{\raise0.3ex\hbox{-}\kern-0.4em D}}

\begin{document}

 \title{Wave equation for generalized Zener model containing complex order fractional derivatives}

 \author{Teodor M. Atanackovi\'c
         \footnote{Institute of Mechanics, Faculty of Technical Sciences, University of Novi Sad,
         Trg D. Obradovi\'ca 6, 21000 Novi Sad, Serbia,
         Electronic mail: atanackovic@uns.ac.rs}\\
         Marko Janev
         \footnote{Institute of Mathematics, Serbian Academy of Sciences and Arts,
         Kneza Mihaila 36, 11000 Belgrade, Serbia,
         Electronic mail: marko\_janev@mi.sanu.ac.rs}\\
         Sanja Konjik
         \footnote{Department of Mathematics and Informatics, Faculty of Sciences, University of Novi Sad,
         Trg D. Obradovi\'ca 4, 21000 Novi Sad, Serbia,
         Electronic mail: sanja.konjik@dmi.uns.ac.rs}\\
         Stevan Pilipovi\'c
         \footnote{Department of Mathematics and Informatics, Faculty of Sciences, University of Novi Sad,
         Trg D. Obradovi\'ca 4, 21000 Novi Sad, Serbia,
         Electronic mail: pilipovic@dmi.uns.ac.rs}
       }

 \date{}
 \maketitle

 \begin{abstract}
 We study waves in a viscoelastic rod whose constitutive equation is of
 generalized Zener type that contains fractional derivatives of complex order.
 The restrictions following from the Second Law of Thermodynamics are
 derived. The initial-boundary value problem for such materials is
 formulated and solution is presented in the form of convolution. Two
 specific examples are analyzed.

 \end{abstract}

 \section{Introduction}
 \label{sec:intro}

Fractional calculus is intensively used for the modelling of various
problems arising in mechanics, physics, engineering, medicine,
economy, biology, chemistry, etc; see  \cite{Podlubny, SKM} 
and references therein as well as our recent books \cite{APSZ-book-vol1, APSZ-book-vol2}. 
Especially, in viscoelasticity, differential operators of arbitrary real order were 
successfully applied, since fractional derivatives, being nonlocal operators, 
 describe intrinsic properties of a material with a "memory", cf.\ \cite{Mainardi10}. 
Derivatives of purely imaginary order were initially studied in \cite{Love-complex},
and later, those of complex order were used to describe viscoelastic 
properties, were studied in \cite{MakrisC-91, MakrisC-92, MakrisC-93}. 
However, in these papers the authors did not consider restrictions for constitutive equations that follow from the 
Second Law of Thermodynamics.
The waves in viscoelastic media have been studied in many papers. We mention few recent studies: In
\cite{Hanyga-16}, wave propagation in viscoelastic bodies, in \cite{Hanyga-07},
nonlinear fractional viscoelastic constitutive equations, while in
 \cite{Tamaogi} the waves have been studied both analytically and experimentally. 
The first fractional generalization of the Zener model for a viscoelastic body was considered in 
\cite{Atanackovic02} while  the distrubuted-order fractional derivatives were considered  
in \cite{APZ10-creep, APZ10-stress}.
Rheological models with fractional damping elements and thermodynamical 
restrictions were studied in \cite{Lion} and those 
with thermal and viscoelastic relaxation effects 
and diffusion phenomena were studied in \cite{Franchi}.
We note that a generalized wave equation, with fractional derivatives, for a viscoelastic material
given in \cite{Wang-16}, is a special case of our
model in this paper.

To the best of our knowledge, so far only factional derivatives of real order 
have been used for describing waves in viscoelastic media. In this work, following
our previous results \cite{AJKPZ15, AKPZ-complexce}, we study waves in a viscoelastic rod whose
material is described by fractional derivatives of complex
order. 
Our particular interest is related to waves in a specific viscoelastic
material described by a generalized Zener standard model. 
A viscoelastic rod model is given by a
system of equations that corresponds to its isothermal motion:
 \begin{align}
 & \frac{\pd}{\pd x}\sigma (x,t) = \rho\,\frac{\pd^{2}}{\pd t^{2}}u(x,t), \notag \\
 & \sigma (x,t) + a_1\, {}_{0}D_{t}^{\alpha}\sigma (x,t) + b_1\, {}_{0}\bar{D}_{t}^{\alpha,\beta}\sigma (x,t)
 = E\,\Big(\eps (x,t) +a_2\, {}_{0}D_{t}^{\alpha}\eps (x,t) + b_2\, {}_{0}\bar{D}_{t}^{\alpha,\beta}\eps (x,t)\Big), \label{1a} \\
 & \eps (x,t) = \frac{\pd}{\pd x}u(x,t), \qquad x\in (0,l],\; \mbox{ or } x\in (0,\infty) \,\, t>0, \notag
 \end{align}
$0<\alpha<1, \beta>0,$ together with the initial conditions
\beq \label{2a}
u(x,0) = 0,\,\,\,\frac{\pd}{\pd t}u(x,0)=0,\,\,\,
\sigma (x,0) = 0,\,\,\,\eps (x,0)=0,
\eeq
and boundary conditions
\beq \label{3a}
u(0,t)=U(t),\,\,\, u(l,t) =0,
\eeq
in the case of finite $l$. In the case $l=\infty$ condition \eqref{3a}$_{2}$
is replaced with
\begin{equation} \label{3b}
\lim_{x\to \infty} u(x,t)=0.
\end{equation}
Here $u,\sigma$ and $\eps$ are displacement, stress and strain,
respectively. Also, $x$ denotes the spatial coordinate oriented
along the axis of the rod and $t$ denotes the time. 

All our calculations are performed on the time variable and $x$ appears as a parameter.
We assume throughout the paper that all the soltions $u$ depending on $x$ and $t$

$$
\mbox{are continuous with respect to } x\in [0,l]  \mbox{ or } x\in [0,\infty), \mbox{ for  almost all } t>0. 
$$
In the case when $u(x,\cdot)$ is a tempered distribution supported by $[0,\infty),$ then we assume tnat
$[0,l]\ni x\mapsto \langle u(x,t),\phi(t)\rangle,$ is continuous for any rapidly decreasing test function $\phi$ ($\phi\in\mathcal S(\R^d)$).

The lenth of the rod is denoted by $l$; in case \eqref{3b}, it is
assumed that $l=\infty$. In the sequel, $[0,l]$ (or $(0,l]$), denotes both cases: with $l<\infty$ and $l=\infty.$

Constants $\rho, E, a_{1}, a_{2}, b_{1},
b_{2}\in \R_{+}$ characterize the material. We note that $E$
represents the modulus of elasticity and $\rho$ density of the
material. The term ${}_{0}D_{t}^{\alpha}$ is the left
Riemann-Liouville fractional derivative operator of order $\alpha$
defined as
$$
{}_{0}D_{t}^{\alpha}u(x,t) = \frac{1}{\Gamma (1-\alpha)}
\frac{d}{dt} \int_{0}^{t} \frac{u(x,\tau)}{(t-\tau)^{\alpha}}d\tau
= \frac{d}{dt} \frac{\tau_+^{-\alpha}}{\Gamma (1-\alpha)}\ast u(x,\tau)(t), 
\quad t>0, \,\,\, 0<\alpha < 1,
$$
where we have assumed that $u(x,\cdot)$ is a polynomially bounded locally integrable function 
supported on $[0,\infty)$, for every $x\geq 0$, and $\ast$ denotes convolution with respect to $t$,
$(f\ast g)(t) =\int_{0}^{t} f(t-\tau) g(\tau)\,d\tau$, $t>0$.

Also in \eqref{1a} we use the following fractional operator of complex order
$$
{}_{0}\bar{D}_{t}^{\alpha,\beta}:=
\frac{1}{2}\Big(\hat{b}_{1}\,{}_{0}D_{t}^{\alpha +i\beta}
+\hat{b}_{2}\,{}_{0}D_{t}^{\alpha -i\beta}\Big),
$$
where the dimensions of constants are chosen to be $\hat{b}_{1}=T^{i\beta}$, $\hat{b}_{2}=T^{-i\beta}$,
so that $|\hat{b}_{1}|=|\hat{b}_{2}|$ ($T$ is a constant having the dimension of time).
This form of a symmetrized fractional derivative of complex order was introduced in
\cite{AJKPZ15, AKPZ-complexce}. Recall that the form of
${}_{0}\bar{D}_{t}^{\alpha,\beta}$ is adapted to the presumption: a fractional derivative of
complex order  applied to a real-valued function has to be 
again real-valued. A dimensionless form is
$$
{}_{0}\bar{D}_{t}^{\alpha,\beta}u(x,t)= \frac{1}{2}\frac{d}{dt} \bigg(
\frac{\tau_+^{-\alpha+i\beta}+\tau_+^{-\alpha-i\beta}}{\Gamma (1-\alpha)}
\bigg) \ast u(x,\tau)(t), 
\quad t>0, \, 0<\alpha < 1, \, \beta>0.
$$

The first equation in \eqref{1a} is the equation of motion
with  $\rho$ being the density of the material. The second
one in \eqref{1a} is the constitutive equation
and coefficients $a_{1},a_{2},b_{1},b_{2}$  satisfy
restrictions that will be determined in the next section. Those restrictions
follow from the Second Law of Thermodynamics. Recall that in the case of the
classical wave equation $u_{tt}=c\,u_{xx}$, for the wave propagation
in an elastic media, the corresponding constitutive equation is given by the
Hooke law $\sigma =E\,\eps$. The last equation in \eqref{1a} is the
strain measure for small local deformations.

Initial conditions \eqref{2a} show that there is no initial displacement, velocity, 
stress and strain, while boundary conditions \eqref{3a}
prescribe displacement at the point $x=0$ and at $x=l$ or infinity.

We introduce dimensionless parameters by a similar consideration as in
\cite{AKOZ-thermodynamical, KOZ10-WZ}. Let
\begin{align*}
& \bar{x} =\frac{x}{L}, \quad
\bar{t}=\frac{t}{T}, \quad
\bar{u}=\frac{u}{L}, \quad
\bar{\sigma}=\frac{\sigma }{E}, \quad
\bar{a}_{i}=\frac{a_{i}}{T^{\alpha}}, \quad
\bar{b}_{i}=\frac{b_{i}}{T^{\alpha}} \,\,\, (i=1,2), \quad
\bar{U}=\frac{U}{L},
\end{align*}
where $T=(a_{2})^{1/\alpha},\ L=(a_{2})^{1/\alpha}\sqrt{\frac{E}{\rho}}$.
Note that $\eps$ is already a dimensionless quantity.
By inserting dimensionless quantities (and dropping the bar sign) into \eqref{1a},
we obtain
\begin{align}
& \frac{\pd}{\pd x}\sigma (x,t)=\frac{\pd^{2}}{\pd t^{2}}u(x,t),  \notag \\
& \sigma (x,t)+a_{1}\,{}_{0}D_{t}^{\alpha}\sigma (x,t)+b_{1}\,{}_{0}\bar{D}_{t}^{\alpha,\beta}\sigma (x,t)
=\eps (x,t)+a_{2}\,{}_{0}D_{t}^{\alpha}\eps(x,t)+b_{2}\,{}_{0}\bar{D}_{t}^{\alpha,\beta}\eps (x,t),  \label{3-sm} \\
& \eps (x,t)=\frac{\pd}{\pd x}u(x,t), \qquad x\in (0,l], \,\,\, t>0. \notag
\end{align}
Assuming that $\sigma(x,\cdot)$ and $\eps(x,\cdot), x\in [0,l)$ are tempered distributions supported by $[0,\infty),$ and applying  the Laplace transform with respect to $t$,
$\cL[\sigma(x,t)](x,s)=\tilde{\sigma}(x,s) = \int_{0}^{\infty} e^{-ts} \sigma(x,t)\,dt$, $\rep s>0$,
(and the same for $\eps$),
to the constitutive equation \eqref{1a}$_{2}$ we obtain
\begin{equation} \label{Lt-1}
\Big(1+a_{1}\,s^{\alpha}+b_{1}\,\big(s^{\alpha +i\beta}+s^{\alpha -i\beta}
\big)\Big)\tilde{\sigma}(x,s)=\Big(1+a_{2}\,s^{\alpha}+b_{2}\,\big(
s^{\alpha +i\beta}+s^{\alpha -i\beta}\big)\Big)\tilde{\eps}(x,s),
\end{equation}
from which we express the stress $\sigma (x,t)$ (after applying the inverse
Laplace transform) as
$$
\sigma (x,t)= \bigg( \cL^{-1}\bigg[\frac{1+a_{2}\,s^{\alpha}+b_{2}\,(s^{\alpha +i\beta}
+s^{\alpha -i\beta})}{1+a_{1}\,s^{\alpha}+b_{1}\,(s^{\alpha +i\beta}
+s^{\alpha -i\beta})}\bigg]\ast_t \eps \bigg) (x,t),
\quad x\in[0,l], \,\,\, t>0.
$$
In the end, this formal calculus will obtain the complete mathematical
justification. 
We use the distributional Laplace transform. Recall that it is defined for 
locally integrable functions of polynomial growth,
and more generally, for tempered distributions supported on $[0,\infty)$.
Important examples are $\cL[\delta(t)](s)=1$ and 
$\cL[H(t)t^\alpha/\Gamma(1+\alpha)](s)=1/s^{\alpha+1}$, $\rep s>0$.
The left hand side has extension for $\rep s=0$ and $|\imp s|\geq \eta_0$,
for any $\eta_0>0$.

If we now replace $\eps$ from \eqref{3-sm}$_{3}$ into \eqref{3-sm}$_{2}$,
and then insert the result into \eqref{3-sm}$_{1}$ we obtain
\begin{equation} \label{cfZwe}
\frac{\pd^{2}}{\pd t^{2}}u(x,t)= L(t) \ast_t
\frac{\pd^{2}}{\pd x^{2}}u(x,t),\qquad x\in [0,l],\,\,\,t>0,
\end{equation}
and the initial and boundary conditions
\begin{eqnarray}
&& u(x,0) =0,\quad \frac{\pd}{\pd t}u(x,0)=0, \qquad x\in[0,l],  \notag \\
&& u(0,t) =U(t),\quad u(l,t) =0 \quad (\mbox{or }\, \lim_{x\to\infty} u(x,t)=0), \qquad t>0,  \label{bc1}
\end{eqnarray}
where
\begin{equation} \label{l-1}
L(t)=\cL^{-1}\bigg[\frac{1+a_{2}\,s^{\alpha}+b_{2}\,(s^{\alpha
+i\beta}+s^{\alpha -i\beta})}{1+a_{1}\,s^{\alpha}+b_{1}\,(s^{\alpha
+i\beta}+s^{\alpha -i\beta})}\bigg](t),
\qquad t>0,\, \alpha \in (0,1),\, \beta >0.
\end{equation}
In the sequel we shall analyze problem \eqref{cfZwe}-\eqref{bc1}. 

Note that it includes several wave equations analyzed earlier. For instance, if the
rod is elastic we have $a_{1}=b_{1}=a_{2}=b_{2}=0$, so that $L(t) =\delta (t)$,
and we obtain the classical wave equation $\frac{\pd^{2}}{\pd t^{2}}u(x,t)
=\frac{\pd^{2}}{\pd x^{2}}u(x,t)$. If the rod is described by a fractional Zener model with
derivatives of real order we have $b_{1}=b_{2}=0$ so that $L(t)=\cL^{-1}
\Big[\frac{1+a_{2}\,s^{\alpha}}{1+a_{1}\,s^{\alpha}}\Big] (t)$, and problem
\eqref{cfZwe}-\eqref{bc1} reduces to the one treated in \cite{KOZ10-WZ} with $l=\infty$.

Restrictions on parameters $a_{1},a_{2},b_{1},b_{2},\alpha,\beta$
will be determined in the following Section, in such a way that the physical meaning
of the problem remains preserved. In this sense, Remark \ref{remv} is important.

\brem
Although it is not obvious, the function $L$ given by \eqref{l-1} is a
real valued function of real variable $t>0$. To see this, we recall the theorem of
Doetsch \cite[p. 293, Satz 2]{Doetsch-Handbuch-1}: A function $L$ is
real-valued (almost everywhere) if its Laplace transform is real-valued for
all real s in the half-plane of convergence on the right from some real $x_{0}$.
Function $\frac{1+a_{2}\,s^{\alpha}+b_{2}\,(s^{\alpha +i\beta
}+s^{\alpha -i\beta})}{1+a_{1}\,s^{\alpha}+b_{1}\,(s^{\alpha +i\beta
}+s^{\alpha -i\beta})}$ clearly satisfies this condition for $s=x>x_{0}=0$
(cf.\ \cite{AKPZ-complexce}).
\erem

The paper is organized as follows:
Following the procedure proposed by Bagley and Torvik 
(see \cite{BagleyTorvik, AKOZ-thermodynamical, AmendolaFabrizioGolden}), 
in the next Section \ref{sec:thermodyn-cfzwe},
we derive thermodynamical restrictions on parameters in \eqref{3-sm}$_{2}$
in order to preserve the Second Law of Thermodynamics. 
Then in Section \ref{sec:estL} we further examine properties of the Laplace transform of
the constitutive equation, which will be needed for the solvability of \eqref{cfZwe}-\eqref{bc1}. 
The existence and uniqueness of a solution to the wave equation \eqref{cfZwe} 
is studied in Section \ref{sec:e-u-w}, where we also explicitly calculate the solution. 
Thermodynamical restrictions again come as essential ones with an appropriate 
sharpness of one of restrictions. 
Results obtained using analytical tools are numerically illustrated in Section \ref{sec:numerics}.

 \section{Thermodynamical restrictions}
 \label{sec:thermodyn-cfzwe}

Consider the constitutive equation \eqref{3-sm}$_2$ for $t>0$, $x\in \R_{+}$,
$\alpha \in (0,1)$ and $\beta >0$. In the analysis that follows the $x$
variable is omitted; \eqref{3-sm}$_2$ is written in the form
\begin{equation} \label{eq1}
\sigma (t)+a_{1}\,{}_{0}D_{t}^{\alpha}\sigma (t)+b_{1}\,{}_{0}\bar{D}_{t}^{\alpha,\beta}\sigma (t)
=\eps (t)+a_{2}\,{}_{0}D_{t}^{\alpha}\eps (t)+b_{2}\,{}_{0}\bar{D}_{t}^{\alpha,\beta}\eps (t).
\end{equation}
We assume that $\sigma(x,\cdot)$ and $\eps(x,\cdot)$ are polynomially bounded locally 
integrable functions supported on $[0,\infty)$, for every $x\geq 0$.
For the coefficients we assume 
$$
a_{i},b_{i}\geq 0, \,\,\, i=1,2, \quad \mbox{ and } \quad a_2>a_1.
$$

Thermodynamical restrictions, i.e., the dissipativity condition - 
the Second Law of Thermodynamics under isothermal conditions -
are closely connected with the following additional assumptions:
 \begin{equation} \label{td-1}
a_{2}b_{1}-a_{1}b_{2}=0,
\end{equation}
\begin{align}
a_{1} &\geq  2b_{1}\cosh \frac{\beta \pi}{2}
\sqrt{1+\Big(\ctg\frac{\alpha \pi}{2}\tgh\frac{\beta \pi}{2}\Big)^{2}}, \notag \\
a_{2} &\geq 2b_{2}\cosh \frac{\beta \pi}{2}
\sqrt{1+\Big(\ctg\frac{\alpha \pi}{2}\tgh\frac{\beta \pi}{2}\Big)^{2}}, \label{td-300}
\end{align}
\begin{align}
a_{1} &\geq  2b_{1}\cosh \frac{\beta \pi}{2}
\sqrt{1+\Big(\tg\frac{\alpha \pi}{2}\tgh\frac{\beta \pi}{2}\Big)^{2}}, \notag \\
a_{2} &\geq 2b_{2}\cosh \frac{\beta \pi}{2}
\sqrt{1+\Big(\tg\frac{\alpha \pi}{2}\tgh\frac{\beta \pi}{2}\Big)^{2}}, \label{td-30}
\end{align}
that will be explained in this section.
We will need a strong inequality in \eqref{td-300} for the existence result in Section 
\ref{sec:estL}, which will be denoted by \eqref{td-300}$_s$.

Thermodynamical restrictions will be determined by following the method proposed 
in \cite{BagleyTorvik}. 

We use the Fourier transform
$$\cF[\vphi(x)](\xi)=\hat{\vphi}(\xi)=\int_{\R} e^{-i\xi x} \vphi(x) dx,\;\xi\in \R \; \mbox{ if } \;\vphi\in L^1(\R),$$ 
or in the sense of tempered distributions if $\vphi$ is locally integrable, supported by
$[0,\infty)$ and bounded by a polynomial. Applying the Fourier transform to \eqref{eq1}
in the sense of distributions (assuming that $\sigma$ and $\eps$ are Fourier transformable),
we obtain $\hat{\sigma}(\omega)=\hat{E}(\omega)\hat{\eps}(\omega)$, so that the
complex modulus of elasticity is
\begin{align*}
\hat{E}(\omega)& =\frac{\hat{P}(\omega)}{\hat{Q}(\omega)}
=\frac{\rep \hat{P}(\omega)+i\imp\hat{P}(\omega)}{\rep \hat{Q}(\omega)+i\imp \hat{Q}(\omega)} \\
& =\frac{\rep \hat{P}(\omega)\rep \hat{Q}(\omega)+\imp \hat{P}(\omega) \imp \hat{Q}(\omega)}{\rep \hat{Q}(\omega)^{2}+
\imp \hat{Q}(\omega)^{2}}
+i\frac{\imp \hat{P}(\omega)\rep \hat{Q}(\omega)-\rep \hat{P}(\omega) \imp \hat{Q}(\omega)}{\rep \hat{Q}(\omega)^{2}+
\imp \hat{Q}(\omega)^{2}} \\
& = \rep\hat{E}(\omega) + i \imp\hat{E}(\omega), \qquad \omega\in (0,\infty), 
\end{align*}
where $\rep\hat{E}(\omega)$ and $\imp\hat{E}(\omega)$ are the loss and the storage modulus. 
Also,
$$
\hat{P}(\omega)=1+a_{2}\,(i\omega)^{\alpha}+b_{2}\,\omega^{\alpha}\Big(
e^{-\frac{\beta \pi}{2}}e^{i(\frac{\alpha \pi}{2}+\ln \omega^{\beta})}
+e^{\frac{\beta \pi}{2}}e^{i(\frac{\alpha \pi}{2}-\ln \omega^{\beta})}
\Big),\quad \omega>0,
$$
$$
\hat{Q}(\omega)=1+a_{1}\,(i\omega)^{\alpha}+b_{1}\,\omega^{\alpha}\Big(
e^{-\frac{\beta \pi}{2}}e^{i(\frac{\alpha \pi}{2}+\ln \omega^{\beta})}
+e^{\frac{\beta \pi}{2}}e^{i(\frac{\alpha \pi}{2}-\ln \omega^{\beta})}
\Big), \quad \omega>0.
$$
With $\hat{P}(0)=\hat{Q}(0)=1,$ $\hat{P}$ and $\hat{Q}$ are continuous function for $\omega\in[0,\infty).$
Let
\begin{eqnarray}
f(\tau,\vphi) &:= & \cos \tau\cos (\alpha \vphi)\cosh (\beta \vphi)
+\sin \tau\sin (\alpha \vphi)\sinh (\beta \vphi),\qquad \tau\in \R, \label{f_def} \\
g(\tau,\vphi) &:= & \cos \tau\sin (\alpha \vphi)\cosh (\beta \vphi)
-\sin \tau\cos (\alpha \vphi)\sinh (\beta \vphi),\qquad \tau\in \R. \label{g_def}
\end{eqnarray}
Properties of functions $f$ and $g$ for $\vphi=\pi/2$ were investigated in \cite{AKPZ-complexce}.
It was shown that the extremal values of $f$ and $g$ are attained at points
$\tau_{f}$ and $\tau_{g}$ respectively, where
$$
\tg \tau_{f}=\tg\frac{\alpha \pi}{2}\tgh\frac{\beta \pi}{2}
\quad \mbox{ and } \quad
\tg \tau_{g}=-\ctg\frac{\alpha \pi}{2}\tgh\frac{\beta \pi}{2}.
$$
There are four such solutions: 
$$
\tau_{f_{1}}\in \Big(0,\frac{\pi}{2}\Big), \tau_{f_{2}}\in \Big(\pi,\frac{3\pi}{2}\Big),
\quad \mbox{ and } \quad
\tau_{g_{1}}\in \Big(\frac{\pi}{2},\pi \Big), \tau_{g_{2}}\in \Big(\frac{3\pi}{2},2\pi \Big).
$$ 
The corresponding extremal values of $f$ and $g$, corresponding to maximum (+)
and minimum (-) are
\begin{align}
f(\tau_{f}, \pi/2) & =\pm \cos \frac{\alpha \pi}{2}\cosh \frac{\beta \pi}{2}
\sqrt{1+\Big(\tg\frac{\alpha \pi}{2}\tgh\frac{\beta \pi}{2}\Big)^{2}}, \notag \\
g(\tau_{g}, \pi/2)& =\pm \sin \frac{\alpha \pi}{2}\cosh \frac{\beta \pi}{2}
\sqrt{1+\Big(\ctg\frac{\alpha \pi}{2}\tgh\frac{\beta \pi}{2}\Big)^{2}}. \label{max-1}
\end{align}

\bprop \label{pr--1}
$\rep\hat{P}(\omega)\geq 1$ and $\rep\hat{Q}(\omega)\geq 1$, $\omega>0$.
\eprop

\pr
We will prove this proposition for $\hat{P}$ since the proof for $\hat{Q}$ 
follows the same lines.
The forms of $P$ and $Q$ imply that their Fourier transform is defined in 
the sense of tempered distributions. Moreover,
$$
\hat{P}(\omega)=\cF[P(t)](\omega)=\tilde{P}(i\omega)=\cL[P(t)](i\omega), 
\quad \omega>0,
$$
where the Laplace transform $\cL[P(t)](z)$ is defined for $\rep z=s\geq 0$.
Next, we have
\begin{align*}
\rep \hat{P}(\omega) &= 1+a_{2}\,\omega^{\alpha} \cos\frac{\alpha \pi}{2}
+ 2b_{2}\,\omega^{\alpha} f(\ln\omega^\beta, \pi/2) \\
& \geq 1+\omega^{\alpha}\cos \frac{\alpha \pi}{2} \bigg(
a_{2}-2b_{2}\cosh \frac{\beta \pi}{2}\sqrt{1+\Big(
\tg\frac{\alpha \pi}{2}\tgh\frac{\beta \pi}{2}\Big)^{2}}
\bigg),
\end{align*}
with $f$ defined by \eqref{f_def} and its minimal value given in \eqref{max-1}$_1$.
Using \eqref{td-30}$_2$ we conclude that $\rep \hat{P}(i\omega)\geq 1$, for all $\omega\in(0,\infty)$.
This proves the proposition.
\ep

The dissipativity condition holds if $\rep\hat{E}(\omega)$ and $\imp\hat{E}(\omega)\geq 0$
for $\omega>0$, see \cite{AmendolaFabrizioGolden, BagleyTorvik}.
These conditions are equivalent to $\rep\hat{P}(\omega)\rep\hat{Q}(\omega)+
\imp\hat{P}(\omega)\imp\hat{Q}(\omega)\geq 0$, and
$\imp\hat{P}(\omega)\rep\hat{Q}(\omega)-\rep\hat{P}(\omega)\imp\hat{Q}(\omega)\geq 0$,
respectively. 
We start with the analysis of $\imp\hat{E}(\omega)\geq 0$, $\omega>0$.
A straightforward calculation yields:
\begin{align}
& \imp\hat{P}(\omega)\rep\hat{Q}(\omega)-\rep\hat{P}(\omega)\imp\hat{Q}(\omega)
=(a_{2}-a_{1})\omega^{\alpha}\sin \frac{\alpha \pi}{2} \notag \\
& \qquad\qquad\qquad
+2(b_{2}-b_{1})\omega^{\alpha}g(\ln \omega^{\beta}, \pi/2)
+2(a_{2}b_{1}-a_{1}b_{2})\omega^{2\alpha}\sin (\ln \omega^{\beta})\sinh \frac{\beta \pi}{2}, \label{g_def1}
\end{align}
where $g$ is given by \eqref{g_def}.

The first observation from \eqref{g_def1} can be stated as follows.

\bprop \label{prop:nstd}
\bi
\item[(i)] A necessary condition for inequality $\imp\hat{E}(\omega)\geq 0$, $\omega >0$, 
is \eqref{td-1}.
\item[(ii)] Necessary and sufficient conditions for both 
$$
\imp\hat{E}(\omega)\geq 0 \quad \mbox{ and } \quad \rep\hat{E}(\omega)\geq 0, 
\quad \omega >0,
$$
are conditions \eqref{td-1}, \eqref{td-300} and \eqref{td-30}.
\ei
\eprop

\pr
(i)\ A careful investigation of \eqref{g_def1} yields that the last term on the right hand side
can take positive and negative values due to the presence of sine function. 
Since it contains the highest power $\omega^{2\alpha}$, while the rest terms in 
\eqref{g_def1} are multiplied by $\omega^{\alpha}$, we conclude that 
$\imp\tilde{E}(\omega)\geq 0$ will be true when $\omega\to\infty$
only if the last term  in \eqref{g_def1} vanishes, i.e., when $a_{2}b_{1}-a_{1}b_{2}=0$.
Thus, this is a necessary condition.

(ii)\ Using \eqref{td-1} in \eqref{g_def1} we obtain
\begin{align*}
& \imp\hat{P}(\omega)\rep\hat{Q}(\omega)-\rep\hat{P}(\omega)\imp\hat{Q}(\omega) \\
& \qquad \geq \omega^{\alpha}\sin \frac{\alpha \pi}{2}\bigg(
(a_{2}-a_{1})-2(b_{2}-b_{1})\cosh \frac{\beta \pi}{2}\sqrt{1+\Big(\ctg\frac{\alpha \pi}{2}\tgh\frac{\beta \pi}{2}\Big)^{2}}
\bigg).
\end{align*}
Thus, $\imp\hat{E}(\omega)\geq 0$, $\omega>0$, if the parameters of the system satisfy
\begin{equation} \label{td-20}
(a_{2}-a_{1})-2(b_{2}-b_{1})\cosh \frac{\beta \pi}{2}\sqrt{1+\Big(
\ctg\frac{\alpha \pi}{2}\tgh\frac{\beta \pi}{2}\Big)^{2}} \geq 0,
\end{equation}
Condition \eqref{td-1} and the assumption $a_2>a_1$ (which together imply $b_2>b_1$) 
when used in \eqref{td-20} leads to \eqref{td-300}. Hence, this is a part of a sufficient 
condition for the dissipativity condition.

Now consider $\rep\hat{E}(\omega)\geq 0$, $\omega>0$.
We calculate:
\begin{align*}
& \rep\hat{P}(\omega)\rep\hat{Q}(\omega)+\imp\hat{P}(\omega)\imp\hat{Q}(\omega) \\
 & \qquad =1+\omega^{\alpha}\Big(
 (a_{1}+a_{2})\cos \frac{\alpha \pi}{2}+2(b_{1}+b_{2})f(\ln \omega^{\beta},\pi/2)
 \Big) \\
& \qquad\quad +\omega^{2\alpha}\Big(
a_{1}a_{2}+2(a_{2}b_{1}+a_{1}b_{2})\cos (\ln\omega^{\beta})\cosh \frac{\beta \pi}{2}
+4b_{1}b_{2}\big(f^{2}(\ln \omega^{\beta},\pi/2)+g^{2}(\ln \omega^{\beta},\pi/2)\big)
\Big),
\end{align*}
where $f$ and $g$ are as in \eqref{f_def} and \eqref{g_def}, respectively.
Since
$$
f^{2}(\ln \omega^{\beta},\pi/2)+g^{2}(\ln \omega^{\beta},\pi/2)
=\cos^{2}(\ln\omega^{\beta})\cosh^{2}\frac{\beta \pi}{2}
+\sin^{2}(\ln \omega^{\beta})\sinh^{2}\frac{\beta \pi}{2}, \quad \omega>0,
$$
we obtain that the third term in the previous equation becomes
\begin{align*}
& a_{1}a_{2}
+2(a_{2}b_{1}+a_{1}b_{2})\cos (\ln \omega^{\beta})\cosh \frac{\beta \pi}{2}
+4b_{1}b_{2}\big(f^{2}(\ln \omega^{\beta},\pi/2)+g^{2}(\ln \omega^{\beta},\pi/2) \\
&\qquad = a_{1}a_{2}
+2a_{2}b_{1}\cos (\ln \omega^{\beta})\cosh \frac{\beta \pi}{2}
+2a_{1}b_{2}\cos (\ln \omega^{\beta})\cosh \frac{\beta \pi}{2} \\
& \qquad \quad +4b_{1}b_{2}\cos^{2}(\ln \omega^{\beta})\cosh^{2}\frac{\beta \pi}{2}
+4b_{1}b_{2}\sin^{2}(\ln \omega^{\beta})\sinh^{2}\frac{\beta \pi}{2} \\
&\qquad=a_{1}a_{2}\bigg(\Big(1+2\frac{b_{1}}{a_{2}}\cos (\ln \omega^{\beta})\cosh \frac{\beta \pi}{2}\Big)^{2}
+4\Big(\frac{b_{1}}{a_{2}}\Big)^{2}\sin^{2}(\ln \omega^{\beta})\sinh^{2}\frac{\beta \pi}{2}\bigg) \\
&\qquad \geq 0, \quad \omega>0,
\end{align*}
where the last equality is derived by the use of condition \eqref{td-1}. Therefore
\begin{align*}
& \rep\hat{P}(\omega)\rep\hat{Q}(\omega)+\imp\hat{P}(\omega)\imp\hat{Q}(\omega) \\
&\qquad \geq
1+\omega^{\alpha}\Big((a_{1}+a_{2})\cos \frac{\alpha \pi}{2}
+2(b_{1}+b_{2})f(\ln \omega^{\beta},\pi/2)\Big) \\
&\qquad \geq 1+\omega^{\alpha}\bigg[ (a_{1}+a_{2})\cos \frac{\alpha \pi}{2}
-2(b_{1}+b_{2})\cos \frac{\alpha \pi}{2}\cosh \frac{\beta \pi}{2}
\sqrt{1+\Big(\tg\frac{\alpha \pi}{2}\tgh\frac{\beta \pi}{2}\Big)^{2}}\bigg], \, \omega>0.
\end{align*}
Thus $\rep\hat{E}(\omega)\geq 0$, $\omega>0$, holds if 
(in fact even more holds, $\rep\hat{E}(\omega)\geq 1>0$)
\begin{equation} \label{td-3}
(a_{1}+a_{2})\geq 2(b_{1}+b_{2})\cosh \frac{\beta \pi}{2}
\sqrt{1+\Big(\tg\frac{\alpha \pi}{2}\tgh\frac{\beta \pi}{2}\Big)^{2}}.
\end{equation}
Again, by inserting condition \eqref{td-1} into \eqref{td-3} and using 
the fact that $a_i,b_i>0$, $a_2,>a_1$, $b_2>b_1$,
we conclude that \eqref{td-3} is equivalent to (\ref{td-30}).

Summing up, we have obtained that inequalities \eqref{td-1}, \eqref{td-300} and \eqref{td-30} 
represent the restrictions following from the Second Law of Thermodynamics, i.e., they are
necessary and sufficient conditions for $\rep\hat{E}(\omega)\geq 0$ and 
$\imp\hat{E}(\omega)\geq 0$ for $\omega \in (0,\infty)$.
\ep

\brem \label{remv}
Consider the assumption $a_i,b_i \geq 0$, $i=1,2$ and \eqref{td-1}.
\begin{enumerate}
\item The case $a_2=a_1$ implies $b_2=b_1$ and this is already explained in Introduction: 
then, $L(t)=\delta(t)$.
\item If $b_1=b_2=0$, then \eqref{td-1} holds for any $a_1, a_2$ and we will not have conditions 
\eqref{td-300} and \eqref{td-30}. Then the needed assumption is $a_2\geq a_1$, cf. \eqref{td-20}.
\item From the beginning we could assume beside \eqref{td-1} and $a_i,b_i\geq 0$, $i=1,2$,  
that $0<a_2<a_1$. This and \eqref{td-1} imply $0<b_2<b_1$. Moreover,  from \eqref{td-1} 
and \eqref{td-20}, we  obtain that 
$$
a_{1} \leq  2b_{1}\cosh \frac{\beta \pi}{2}
\sqrt{1+\Big(\ctg\frac{\alpha \pi}{2}\tgh\frac{\beta \pi}{2}\Big)^{2}}, 
$$
\begin{equation} \label{td-400}
a_{2} \leq 2b_{2}\cosh \frac{\beta \pi}{2}
\sqrt{1+\Big(\ctg\frac{\alpha \pi}{2}\tgh\frac{\beta \pi}{2}\Big)^{2}}. 
\end{equation}
So, in this case, \eqref{td-1}, \eqref{td-400} and \eqref{td-30}
are necessary and sufficient thermodynamical conditions.
\end{enumerate}
\erem

 \section{Necessary estimates for the Laplace transform}
 \label{sec:estL}

Note that the complex modulus $\hat{E}(\omega)$ may be treated as a special
value of the complex function $\tilde{E}(s)$ obtained from the Laplace transform 
of the constitutive equation (now written with both arguments) 
$\tilde{\sigma}(x,s)=\tilde{E}(s)\tilde{\eps}(x,s)$, $\rep s \geq 0$,
see \eqref{Lt-1}, calculated for $\rep s=0$, $\imp s=\omega \in \R_{+}$.

In the sequel we examine properties of $\tilde{E}(s)$. Actually, we will work with $M^2=1/\tilde{E}$.
Moreover, in the rest of the paper we will assume $a_2>a_1>0$, $b_2>b_1>0$,
\eqref{td-1}, \eqref{td-30} and the stronger assumption \eqref{td-300}$_s$. 
This will be an essential point of the second part of the proof of Proposition \ref{prop-s}.

From \eqref{Lt-1} we have, formally,
\beq \label{so1a}
M^2(s) = \frac{1}{\tilde{E}(s)}
= \frac{1+a_{1}\,s^{\alpha}+b_{1}\,\big[s^{\alpha +i\beta}+s^{\alpha -i\beta}
\big]}{1+a_{2}\,s^{\alpha}+b_2\,\big[s^{\alpha +i\beta}+s^{\alpha -i\beta}\big]}
= \frac{\tilde{Q}(s)}{\tilde{P}(s)}, \qquad \rep s >0.
\eeq
With $M^2(0)=1$, $M^2$ and its square root $M$ are continuous functions on $\rep \geq 0.$
We will show in Proposition \ref{prop:noMsing} that $\tilde{P}(s)$ does not have 
zeros in the domain $\rep s\geq 0$. This means that $M^2$ is well defined.
We will use
$$
\rep M^2=\rep\tilde{E}/|\tilde{E}|^2, \qquad
\imp M^2=-\imp\tilde{E}/|\tilde{E}|^2.
$$
So calculating $\tilde{E}$ we will get estimates of $M^2$.
Let $s=\rho e^{i\vphi}=s_{0}+ip$, with $\rho =\sqrt{s_{0}^{2}+p^{2}}$, $s_{0},p\geq 0$, so
that $\vphi \in [0,\pi/2]$. Then, we compute
\begin{eqnarray}
\rep\hat{E}(s) &=&\frac{1+B+C}{D}  \notag \\
\imp\hat{E}(s) &=&\frac{\rho^{\alpha}\Big(
(a_{2}-a_{1})\sin (\alpha\vphi)+2(b_{2}-b_{1})g(\ln \rho^{\beta},\vphi)
\Big)}{D}, \label{Lap-1}
\end{eqnarray}
where ($\rep s=s_{0}$, $p\geq 0$, $\vphi \in [0,\pi/2]$)
\begin{align*}
B &= \rho^{\alpha}\Big(
(a_{1}+a_{2})\cos (\alpha \vphi)+2(b_{1}+b_{2})f(\ln \rho^{\beta},\vphi)
\Big), \\
C &= \rho^{2\alpha}\Big[
a_{1}a_{2}+4a_{1}b_{2}\cos (\ln \rho^{\beta})\cosh (\beta \vphi)
+ 4b_{1}b_{2}(f^{2}(\ln \rho^{\beta},\vphi)+g^{2}(\ln \rho^{\beta},\vphi)
\Big], \\
D &=1+ 2\rho^{\alpha}(a_{1}\cos (\alpha \vphi)+2b_{1}f(\ln \rho^{\beta},\vphi)) 
+ \rho^{2\alpha}\Big[
a_{1}^{2}+4a_{1}b_{1}(\cos (\alpha \vphi)f(\ln \rho^{\beta},\vphi) \\
& \qquad
+\sin (\alpha \vphi)g(\ln \rho^{\beta},\vphi))+4b_{1}^{2}(f^{2}(\ln \rho^{\beta},\vphi))+g^{2}(\ln \rho^{\beta},\vphi)
\Big],
\end{align*}
with $f$ and $g$ as in \eqref{f_def} and \eqref{g_def}, respectively.

We have, in the domain $\rep s=s_{0}$, $p\geq 0$, that \eqref{Lap-1} implies
\begin{equation} \label{est-1}
\imp M^2 \sim c\rho^{-2\alpha}, \,\,\, c>0, \,\,\, \mbox{ as } \,\,\, \rho\to \infty.
\end{equation}
 
Following the same procedure as in \cite{AKPZ-complexce}
in obtaining \eqref{max-1}, we conclude that the
maximal and minimal values of $f(\tau,\vphi)$ and $g(\tau,\vphi)$ with respect
to $\tau$ are
\begin{align}
f(\tau_{f},\vphi)& =\pm \cos (\alpha \vphi) \cosh (\beta \vphi) \sqrt{1+\big(
\tg (\alpha \vphi) \tgh (\beta \vphi) \big)^{2}}, \notag \\
g(\tau_{g},\vphi)& =\pm \sin (\alpha \vphi) \cosh (\beta \vphi) \sqrt{1+\big(
\ctg (\alpha \vphi) \tgh (\beta \vphi) \big)^{2}}.  \label{Lap-3}
\end{align}


\bprop \label{prop-s}
Suppose that the thermodynamical restrictions 
\eqref{td-1}, \eqref{td-300}$_s$ and \eqref{td-30}
are satisfied (as well as conditions on $a_i, b_i$, $i=1,2$). 
Let $s_{0}\geq 0$ be fixed and $s=s_{0} + ip$, $p\in\R$.
Then $\rep M^2(s) >0$, for all $ p>0$. Also, there exists $p_{0}>0$ such that 
$\imp M^2(s)<0$, $s=s_{0} + ip$, $p>p_{0}$.
\eprop

\pr 
Set $s=\rho e^{i\vphi}$, $\rho=\sqrt{s_{0}^2+p^2}$, $\vphi\in[0,\pi/2]$.
By \eqref{f_def} and \eqref{g_def} we have that $\rep\tilde{E}(s)>0$ if 
$1+B+C> 0$. Using $f^{2}(\ln \rho ^{\beta },\varphi )+g^{2}(\ln \rho ^{\beta
},\varphi )=\cos ^{2}(\ln \rho ^{\beta })\cosh ^{2}\beta \varphi
+\sin ^{2}(\ln \rho ^{\beta })\sinh ^{2}\beta \varphi $, we obtain
\begin{align*}
C & = \rho^{2\alpha}\Big[ 
a_{1}a_{2}+4a_{1}b_{2}\cos (\ln \rho^{\beta})\cosh(\beta \vphi) 
+4b_{1}b_{2}(f^{2}(\ln \rho^{\beta},\vphi)+g^{2}(\ln \rho^{\beta},\vphi))
\Big] \\
& = \rho^{2\alpha}\bigg[ 
a_{1}a_{2}\bigg(\Big(1+2\frac{b_{2}}{a_{2}}\cos (\ln\rho^{\beta})\cosh \beta \vphi \Big)^{2}
+4\Big(\frac{b_{1}}{a_{2}}\Big)^{2}\sin^{2}(\ln \rho^{\beta})\sinh^{2}\beta \vphi \bigg)\bigg] \\
& \geq 0.
\end{align*}
To estimate $B$ we use \eqref{Lap-3}$_{1}$ so that
\begin{align*}
&(a_{1}+a_{2})\cos (\alpha \vphi)+2(b_{1}+b_{2})f(\ln \rho^{\beta},\vphi) \\
&\qquad \geq \cos (\alpha \vphi) \Big( 
(a_{1}+a_{2})-2(b_{1}+b_{2})\cosh (\beta\vphi)\sqrt{1+\tg^{2}(\alpha \vphi) \tgh^{2}(\beta \vphi)}
\Big).
\end{align*}
Since $\cosh (\beta \vphi)\sqrt{1+\tg^{2}(\alpha\vphi)\tgh^{2}(\beta \vphi)}
<\cosh \frac{\beta\pi}{2}\sqrt{1+\tg^{2}\frac{\alpha \pi}{2}\tgh^{2}\frac{\beta \pi}{2}}$
(all functions $\cosh$, $\tg$ and $\tgh$ are monotone increasing functions for $\vphi\in[0,\pi/2]$),
\eqref{td-3} with \eqref{td-30} implies that $B\geq 0$. Therefore, $\rep\tilde{E}(s)> 0$.

To estimate $\imp\tilde{E}(s)$ we start from \eqref{Lap-1}$_{2}$ and
analyze the term $(a_{2}-a_{1})\sin (\alpha \vphi)+2(b_{2}-b_{1})g(\ln \rho^{\beta},\vphi)$.
Note that
\begin{align*}
& (a_{2}-a_{1})\sin (\alpha \vphi)+2(b_{2}-b_{1})g(\ln \rho^{\beta},\vphi) \\
&\qquad \geq (a_{2}-a_{1})\sin (\alpha \vphi)
+2(b_{2}-b_{1}) \min_{x\in R} g(x,\vphi) \\
&\qquad = \sin (\alpha \vphi)\Big[ 
(a_{2}-a_{1})-2(b_{2}-b_{1})\cosh (\beta \vphi) \sqrt{1+\big(\ctg (\alpha \vphi) \tgh (\beta \vphi) \big)^{2}}
\Big].
\end{align*}
Thus $\imp\tilde{E}(s)> 0$ if
\begin{equation} \label{est-3}
(a_{2}-a_{1})> 2(b_{2}-b_{1})\cosh (\beta \vphi) \sqrt{1+\big(\ctg (\alpha \vphi) \tgh (\beta \vphi) \big)^{2}}.
\end{equation}
Relation \eqref{est-3} becomes \eqref{td-20} with a strict inequality when $\vphi = \pi/2$. 
Since $\vphi \to \pi/2$ when $p\to \infty$, and right hand side of \eqref{est-3} is
continuous function of $\vphi$, we conclude that there is $p_{0}> 0$ such that 
\eqref{est-3} is satisfied for all $s=s_{0}+ip$ with $p>p_{0}$. 
Therefore, $\imp\tilde{E}(s)> 0$, $s=s_{0}+ip$, $p>p_{0}$. 
\ep

\brem \label{rem:cconj} 
Using the symmetry properties of trigonometric and hyperbolic functions, and 
the fact that $f(x,-\vphi)=f(x,\vphi)$ and $g(x,-\vphi)=-g(x,\vphi)$, we conclude
that $\tilde{E}(\bar{s})=\overline{\tilde{E}(s)}$. Thus, $\rep \tilde{E}(s)>0$, 
for all $s=s_{0}-ip$ satisfying $s_0\geq 0$, $p>0$. Also, $\imp \tilde{E}(s)< 0$, 
for all $s=s_{0}-ip$ satisfying $s_0\geq 0$, $p>p_0$. 
\erem

 \section{Solution to \eqref{cfZwe}-\eqref{bc1}}
 \label{sec:e-u-w}

Recall that we continue to use conditions \eqref{td-1}, \eqref{td-300}$_s$ and \eqref{td-30}$_s$.

We return to the initial-boundary value problem \eqref{cfZwe}-\eqref{bc1}. 
Applying  the Laplace transform to \eqref{cfZwe} we obtain
\begin{equation} \label{so1}
\frac{d^{2}\hat{u}(x,s)}{dx^{2}}-s^{2}M^{2}(s) \hat{u}(x,s)=0, \quad x>0, \; \rep s>0,
\end{equation}
where $M^{2}(s)$ is defined by \eqref{so1a}.
Boundary conditions become
\begin{equation} \label{so2}
\tilde{u}(0,s)=\tilde{U}(s), \quad \tilde{u}(l,s)=0, \mbox{ if } l \mbox{ is finite},
\end{equation}
\begin{equation} \label{so3}
\tilde{u}(0,s)=\tilde{U}(s), \quad \lim_{x\to\infty} \tilde{u}(x,s)=0, \mbox{ if } l=\infty.
\end{equation}
Solutions to \eqref{so1}, \eqref{so2} and \eqref{so1}, \eqref{so3} are
\beq \label{so4}
\tilde{u}(x,s) = \tilde{U}(s) \bigg[
\frac{e^{sM(s) x}}{1-e^{2sM(s) l}}
+
\frac{e^{-sM(s) x}}{1-e^{-2sM(s) l}}
\bigg], x>0, \quad \rep s>0,
\eeq
and
\begin{equation} \label{so5}
\tilde{u}(x,s)=\tilde{U}(s) e^{-sM(s)x}, \quad x>0,\; \rep s>0,
\end{equation}
respectively.

We need the following result on $M$.

\bprop \label{prop:noMsing}
$M(s)$ has no singular points with positive real part.
\eprop

\pr
Singular points of $M(s)$ are zeros of $\tilde{P}(s)$, $\rep s>0$, see \eqref{so1a}.
Thus, we consider the equation
$\tilde{P}(s) \equiv 1+a_{2}\,s^{\alpha}+b_{2}\,\big[ s^{\alpha+i\beta}+s^{\alpha -i\beta}\big] =0$.
To determine zeros of $\tilde{P}$ we use the argument principle.
Note that if $s_{0}$ is a solution to $\tilde{P}(s) =0$ then the complex conjugate
$\bar{s}_{0}$ is also a solution since $\tilde{P}(\bar{s})
=1+a_{2}\,\bar{s}^{\alpha}+b_{2}\,\big[\bar{s}^{\alpha +i\beta}
+\bar{s}^{\alpha -i\beta}\big] =\overline{\tilde{P}(s)}$ (see also Remark \ref{rem:cconj}).
Therefore it is enough to consider zeros in the part of the complex plane with
$\rep s\geq 0$, $\imp s\geq 0$.

 \begin{figure}[htb]
 \begin{center}
 \includegraphics[width=5.cm]{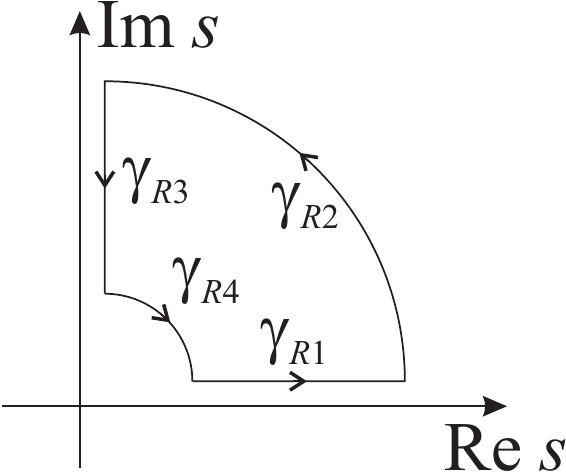}
 \end{center}
 \caption{Integration path $\Gamma$ for zeros of $\hat{P}$}
 \label{fig:kontura1}
 \end{figure}

Let $\Gamma =\gamma_{R1}\cup \gamma_{R2}\cup
\gamma_{R3}\cup \gamma_{R4}$ be a contour as shown in Figure \ref{fig:kontura1}.

Contour $\gamma_{R1}$ is parametrized by $s=x$, $x\in(\eps,R)$ with $\eps\to 0$,
$R\to \infty$. On this part we have $\imp \tilde{P}(s) =0$, $\rep \tilde{P}(s)
=1+x^{\alpha}(a_{2}+2b_{2}\cos (\ln x^{\beta}) )\geq 1+x^{\alpha}(a_{2}-2b_{2})$.
Thus, $\rep \tilde{P}(s) \to \infty$ as $x\to \infty$ if $a_{2}>2b_{2}$,
which is satisfied due to the thermodynamical restrictions \eqref{td-30}$_{2}$,
see the end of Proposition \ref{prop:nstd}.

Along $\gamma_{R2}$ we have $s=Re^{i\vphi}$, $\vphi \in [0,\pi/2]$, $R\to \infty$,
so that
\begin{align*}
\rep \tilde{P}(s) &= 1+a_{2}R^{\alpha}\cos (\alpha\vphi)
+2b_{2}R^{\alpha} f(\ln R^\beta, \vphi)
\\
\imp \tilde{P}(s) &= a_{2}R^{\alpha}\sin (\alpha \vphi)
+2b_{2}R^{\alpha} g(\ln R^\beta, \vphi),
\end{align*}
where $f$ and $g$ are given by \eqref{f_def} and \eqref{g_def}.
The minimum of $f$ is given by \eqref{Lap-3}, hence $\rep \tilde{P}(s)> 1$
on $\gamma_{R2}$ if
$$
a_{2}> 2b_{2}\cosh \frac{\beta \pi}{2}\sqrt{1+\Big(\tg \frac{\alpha\pi}{2}\tgh\frac{\beta \pi}{2}\Big)^{2}}.
$$
This condition is satisfied because of the dissipation inequality \eqref{td-300}$_s$.
Moreover, we have
\begin{align*}
\rep \tilde{P}(s) &\to \infty \,\,\, \mbox{ and } \,\,\, \imp \tilde{P}(s) \to 0,
\,\,\, \mbox{ for } \,\,\, \vphi =0,\, R\to \infty, \\
\rep \tilde{P}(s) &\to \infty \,\,\, \mbox{ and } \,\,\, \imp \tilde{P}(s) \to \infty,
\,\,\, \mbox{ for } \,\,\,  \vphi =\frac{\pi}{2}, \, \alpha <1, \, R\to \infty.
\end{align*}

On $\gamma_{R3}$ we have $s=ip$, $p\in (\eps,R)$ with $\eps \to 0$, $R\to \infty$.  
According to the calculated Fourier transform 
of $P$ and its relation with the Laplace transform, we have
\begin{align*}
\rep \tilde{P}(s=ip) &=1+a_{2}p^{\alpha}\cos \frac{\alpha \pi}{2}
+ 2b_{2}p^{\alpha}f(p,\pi/2) \\
&\geq 1+p^{\alpha}\cos \frac{\alpha \pi}{2} \bigg(
a_{2}-2b_{2}\cosh \frac{\beta \pi}{2}\sqrt{1+\Big(
\tg\frac{\alpha \pi}{2}\tgh\frac{\beta \pi}{2}\Big)^{2}}
\bigg),
\end{align*}
with $f$ defined by \eqref{f_def} and its minimal value given in \eqref{max-1}$_1$.
Using \eqref{td-30}$_2$ we conclude that $\rep \tilde{P}(s=ip)\geq 1$, for all $p\in (\eps,R)$.

Finally, parametrization of $\gamma_{R4}$ is $s=\eps e^{i\vphi}$,
$\vphi \in [0,\pi/2]$, $\eps \to 0$, so that $\rep \tilde{P}(\eps e^{i\vphi}) \to1$,
$\imp \tilde{P}(\eps e^{i\vphi}) \to 0$, as $\eps\to 0$.

All together, the change of the argument of function $\tilde{P}$ along
$\Gamma$ is zero, i.e.,
$$
\Delta \arg \psi (s) =0,
$$
implying that there are no zeroes of $\tilde{P}$ in the right complex half-plane.
Therefore there are no singular points of $M(s)$ with $\rep s>0$.
\ep

The following Corollary provides that the solution \eqref{so4} is well defined.

\bcor
Functions  $1-exp (2sM(s)l)$ and $1-exp (-2sM(s)l)$ have no zeros in the right
complex half-plane $\rep s\geq 0$.
\ecor

\pr
Zeros of functions $1-exp (2sM(s)l)$ and $1-exp (-2sM(s)l)$ satisfy $sM(s)=0$.
The claim now follows from Proposition \ref{prop:noMsing}.
\ep

\brem \label{rem:signM} 
Using our calculations for the Fourier transformation of 
$E$, we have
$$
M^{2}(i\omega) = \frac{1}{\hat{E}(\omega)}
=\frac{\rep\hat{E}(\omega)-i\imp\hat{E}(\omega)}{|\hat{E}(\omega)|^{2}},
$$
and $\rep\hat{E}(\omega), \imp\hat{E}(\omega)\geq 0$ for $\omega >0$, we
conclude that $\rep M^{2}(i\omega) >0$ and $\imp M^{2}(i\omega)< 0$ for 
$\omega>0$ if the thermodynamical restrictions \eqref{td-1}, 
\eqref{td-300}$_s$ and \eqref{td-30} are satisfied.
Hence $\rep M(i\omega) > 0$ and $\imp M(i\omega) < 0$, $\omega>0$. Also,
from Remark \ref{rem:cconj}, it follows $\rep M(-i\omega) > 0$ and
$\imp M(-i\omega) < 0$, $\omega>0$.

Similarly, from 
$$
M^{2}(s) =\frac{1}{\hat{E}(s)}=\frac{\rep\hat{E}(s)-i\hat{E}(s)}{|\hat{E}(s)|^{2}}
$$ 
and Proposition \ref{prop-s} we conclude that $\rep M^2(s)>0$, $s = s_0+ip$, 
$s_0,p>0$, and $\imp M^{2}(s) \leq 0$, $s = s_0+ip$, $s_0>0$, $p>p_0$.
Consequently, $\rep M(s)> 0$, $s = s_0+ip$, $s_0,p>0$, and 
$\imp M(s) < 0$, $s = s_0+ip$, $s_0>0$, $p>p_0$. 
\erem

We now state the main result of this Section.

\bthm 
Problem \eqref{cfZwe}-\eqref{bc1} has a solution given as
\begin{equation} \label{res-0}
u(x,t)= (U\ast_t K) (x,t) =\int_{0}^{t} U(t-\tau) K(x,\tau) \,d\tau, \quad x\in (0,l], \,\,\, t>0,
\end{equation}
where, for $l=\infty,$
$$
K(x,t) = \frac{1}{2\pi i} \int_{s_{0}-i\infty }^{s_{0}+i\infty}
\exp (ts) \Big[
\frac{\exp (sM(s)x)}{1-exp (2sM(s)l)}
+ \frac{\exp (-sM(s)x)}{1-exp (-2sM(s)l)} \Big] \,ds,
$$
or
\begin{equation} \label{res-1a}
K(x,t) =\frac{1}{2\pi i} \int_{s_{0}-i\infty }^{s_{0}+i\infty} \exp (ts) \exp (-sM(s)x) \,ds
\end{equation}
if $l=\infty$. Here $s_{0}>0$.
In particular, for every $t$, $K$ is continuous with respect to $x$, as well as $u$.
If $x=0$ then $K(0,t)=\delta(t)$.
\ethm

In the case $l=\infty$ solution \eqref{res-0} can be computed more explicitly,
as it is given in the following statement.

\bthm \label{th:sol-infty}
Let $l=\infty$. Then the solution kernel \eqref{res-1a} takes the form
\begin{equation} \label{res-10}
K(x,t)=\frac{1}{\pi} \int_{0}^{\infty} \exp \big[ \tau \imp M(i\tau) x \big]
\cos \big[ \tau (t-\rep M(i\tau) x) \big] \,d\tau, \quad x> 0, \,\,\, t>0.
\end{equation}
It is continuous function of $x$. If $x=0$ then $K(0,t)=\delta(t)$.
\ethm

\pr
Set $s=s_{0}+i\tau$ ($ds=i\,d\tau$) and $p>0$. Then \eqref{res-1a} becomes
$$
K(x,t)=\frac{1}{2\pi i} \lim_{p\to \infty} \int_{-p}^{p}
\exp \big[(s_{0}+i\tau) (t-M(s_{0}+i\tau) x)\big] \,d\tau, \quad x> 0, \,\,\, t>0.
$$
The estimate below shows the continuity with respect to $x>0$ for every $t>0$.

\begin{figure}[htb]
 \begin{center}
 \includegraphics[width=5.cm]{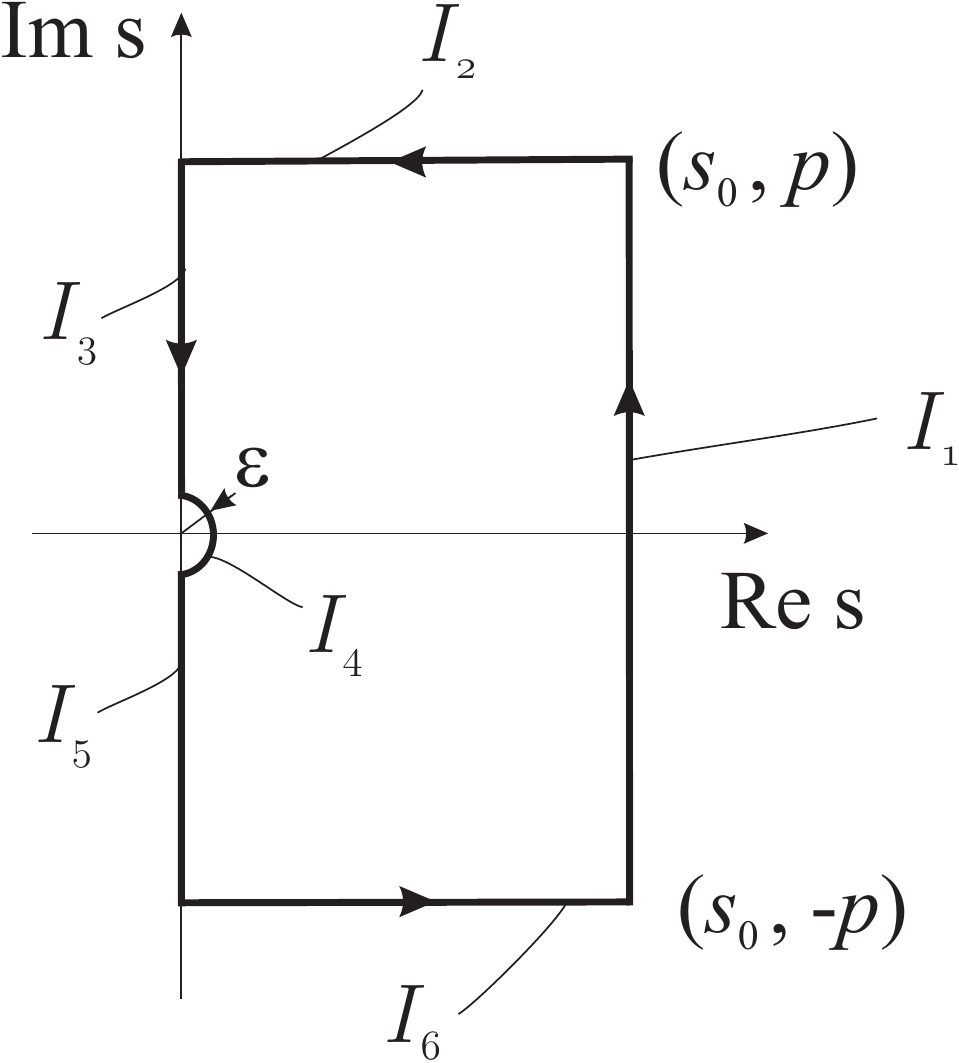}
 \end{center}
 \caption{Integration path $I$}
 \label{fig:kontura2}
 \end{figure}

Consider the contour shown in Figure \ref{fig:kontura2}. Then
$K(x,t) =\frac{1}{2\pi} \lim_{p\to \infty} I_{1} (x,t,p)$, $x\geq 0$, $t>0$, where
the Cauchy integral theorem implies
$$
I_{1}=-(I_{2}+I_{3}+I_{4}+I_{5}+I_{6}).
$$

For the integral $I_{2}$ we have to use \eqref{est-1} and 
\begin{align*}
|I_{2}| &= \Big| -\int_{0}^{s_{0}} \exp \big[(\sigma +ip) (t-M(\sigma +ip)x)\big] \,d\sigma \Big| \\
& \leq \int_{0}^{s_{0}} \exp \big[\sigma(t-\rep M(\sigma+ip)x)+p\imp M(\sigma +ip) x\big] \,d\sigma \\
& \leq C\, \int_{0}^{s_{0}} \exp \big[p\imp M(\sigma +ip)\big] \,d\sigma \\
& \leq C\, \int_{0}^{s_{0}} \exp \big[-p^{1+\frac{\alpha}{2}}\big] \,d\sigma \\
& < \infty
\end{align*}

Similar arguments prove that $\lim_{p\to\infty} I_{6}=0$.

Next for $I_{4}$ we have $s=\eps \exp (i\vphi)$, $ds=i\eps \exp (i\vphi) \,d\vphi$
so that
$$
\lim_{\eps \to 0} I_{4}=\lim_{\eps\to 0} \int_{-\pi/2}^{\pi/2} \exp \big[\eps \exp(i\vphi)
(t-M(\eps \exp (i\vphi)) x) \big] i\eps \exp (i\vphi) \,d\vphi =0.
$$
Therefore $I_{1}=- (I_{3}+I_{5})$ so that with $s=i\tau$ we get
\begin{align*}
I_{1} &= -i\Big[ \int_{p}^{\eps} \exp \big[ i\tau (t-M(i\tau)x) \big] \,d\tau
+ \int_{-\eps }^{-p} \exp \big[ i\tau (t-M(i\tau)x) \big] \,d\tau \Big] \\
&= -2i\Big[ \int_{p}^{\eps} \exp \big[ \tau \imp M(i\tau)x \big]
\cos \big[ \tau (t-\rep M(i\tau)x) \big] \,d\tau \Big],
\end{align*}
where we used $\imp M(-ip) =-\imp M(ip)$, $\rep M(ip) =\rep M(-ip)$.
Thus,
$$
K(x,t)=\frac{1}{\pi} \lim_{\eps \to 0} \Big(
\lim_{p\to \infty} \int_{\eps }^{p} \exp \big[ \tau \imp M(i\tau) x \big]
\cos \big[ \tau (t-\rep M(i\tau) x) \big] \,d\tau \Big), \quad x> 0, \,\,\, t>0,
$$
which proves the claim.
\ep

\brem 
In the numerical example of the last Section we can notify an 
oscillatory character of the solution which comes from integrals $I_2$ and $I_6$  
that is obtained in numerical experiments for relatively small values of 
$p$ (see the next Section \ref{sec:numerics}). 
\erem

 \section{Numerical experiments}
 \label{sec:numerics}

\begin{figure}[htb]
 \begin{center}
 \includegraphics[width=6.cm]{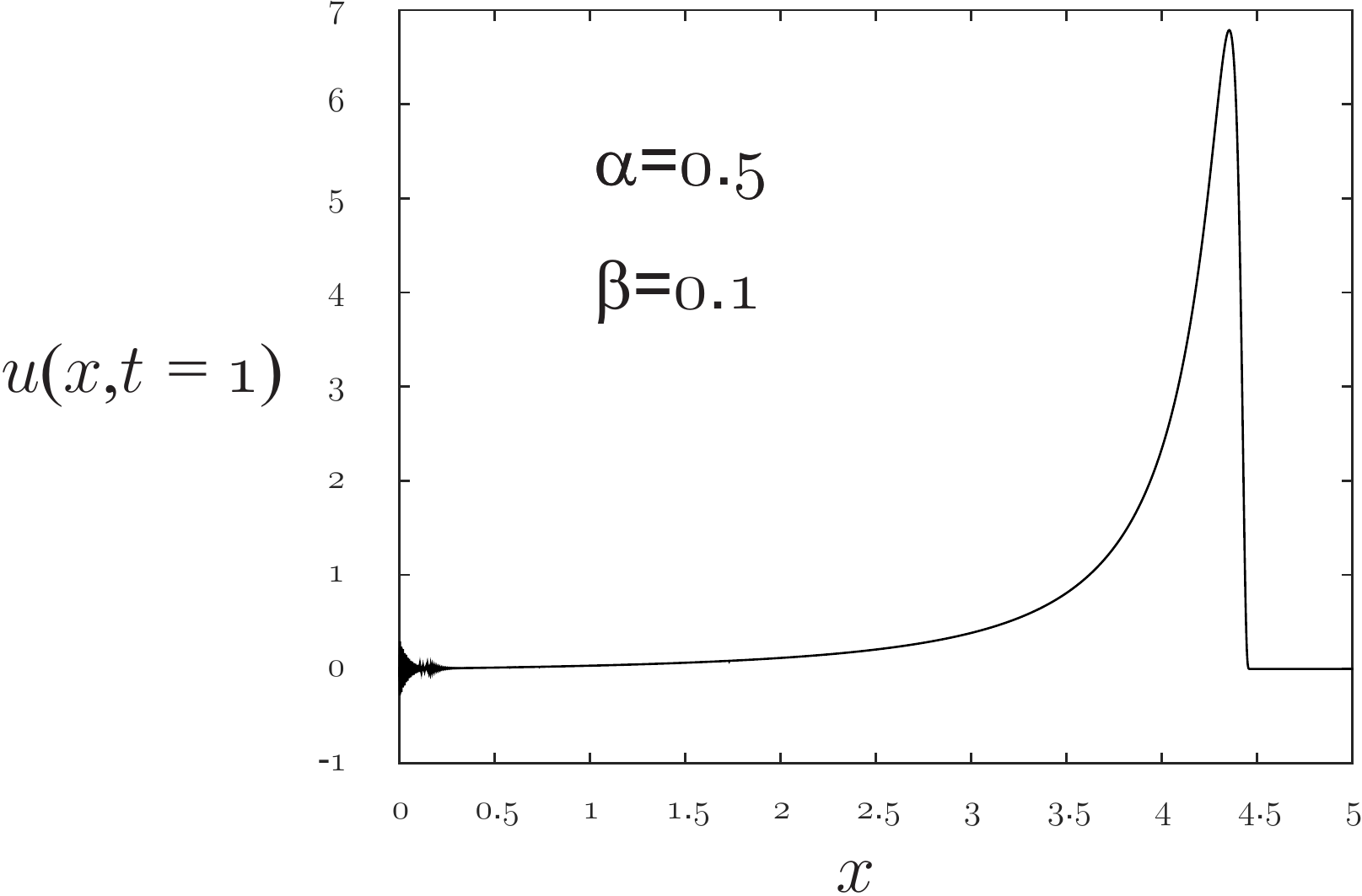}
 \end{center}
 \caption{Displacement $u$ for $\alpha =0.5$, $\beta =0.1$ and $t=1$.}
 \label{fig:primer1}
 \end{figure}

Results obtained in previous sections will be presented for various sets of 
parameters and boundary conditions. The goal is to examine and compare 
how different values of coefficients and orders of fractional derivatives
influence the solution, i.e., how the solution wave propagates in the viscoelastic 
rod. We treat two different cases:

Case 1.  Choose the following values for the parameters in system \eqref{cfZwe}-\eqref{bc1}
$$
a_1=1, \quad
a_2=20, \quad
b_1=0.1, \quad
l=\infty.
$$
Then according to \eqref{td-1}, $b_2=a_2b_1/a_1=2$, and one checks easily, by inserting these chosen
values of parameters into \eqref{td-300} and \eqref{td-30}, that the thermodynamical restrictions
are satisfied.
Suppose that $U(t) =\delta (t)$, the Dirac distribution. Then combining \eqref{res-0} and
\eqref{res-10} the solution reads
\begin{equation} \label{res-11}
u(x,t)=\frac{1}{\pi} \int_{0}^{\infty} \exp \big[ \tau \imp M(i\tau)x \big]
\cos \big[ \tau ( t-\rep M(i\tau)x) \big] \,d\tau.
\end{equation}
In Figure \ref{fig:primer1} we present solution $u$ given by \eqref{res-11} 
for $\alpha=0.5$, $\beta = 0.1$ and $t= 1$.

\begin{figure}[htb]
 \begin{center}
 \includegraphics[width=6.cm]{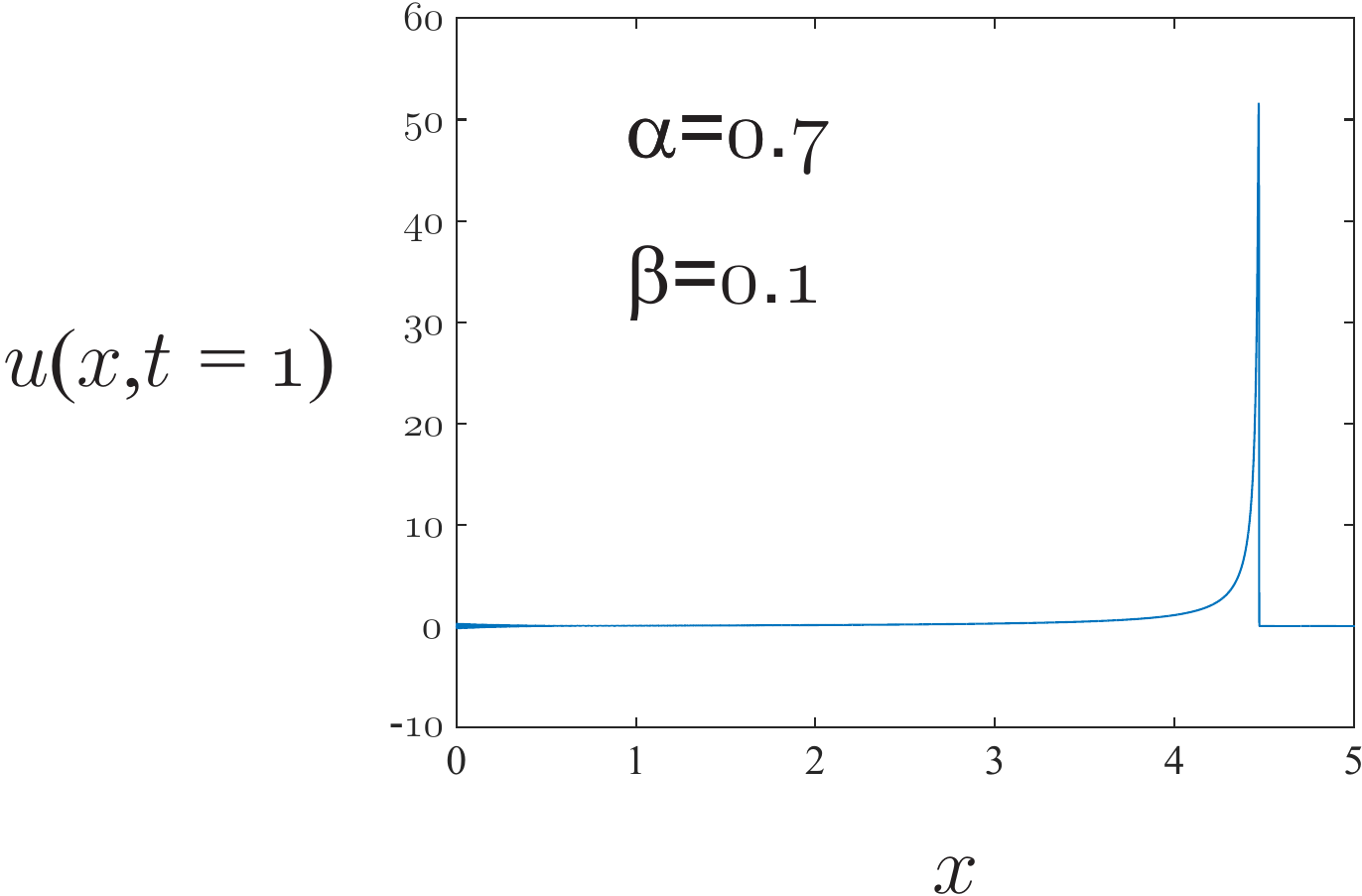}
 \end{center}
 \caption{Displacement $u$ for $\alpha =0.7$, $\beta =0.1$ and $t=1$.}
 \label{fig:primer2}
 \end{figure}

In order to examine the influence of $\alpha$ on the solution, in Figure \ref{fig:primer2} 
we present $u$ given by \eqref{res-11} for $\alpha=0.7$, $\beta = 0.1$ and $t= 1$.

\begin{figure}[htb]
 \begin{center}
 \includegraphics[width=6.cm]{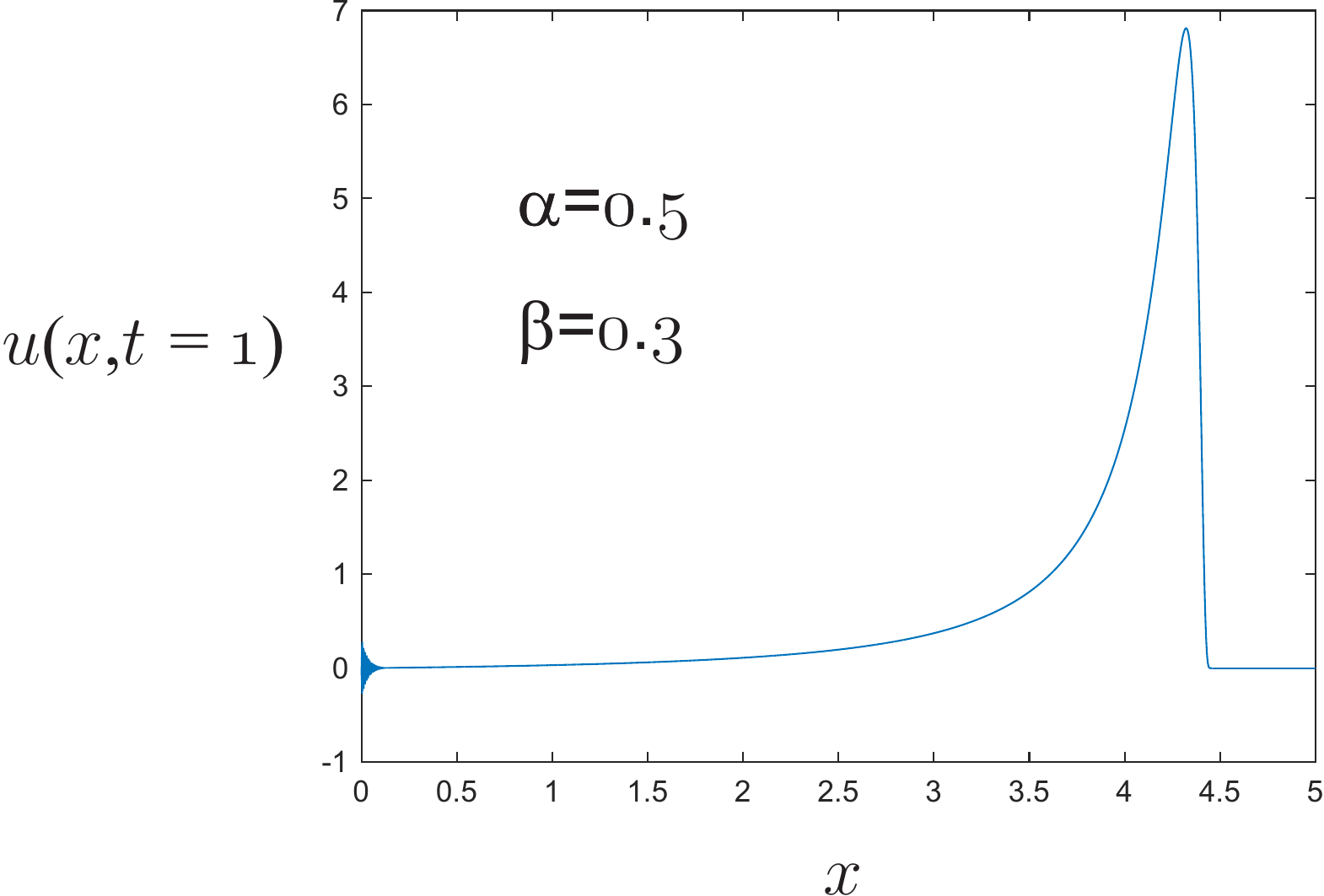}
 \end{center}
 \caption{Displacement $u$ for $\alpha =0.5$, $\beta =0.3$ and $t=1$.}
 \label{fig:primer3}
 \end{figure}

Finally, we increase $\beta$, so we consider the case $\alpha =0.5$, $\beta = 0.3$ and $t= 1$. 
The results are shown in Figure \ref{fig:primer3}.

Case 2. Suppose now that $U(t) =H(t)$, where $H$ is the Heaviside function. 
Also, we assume that the parameters are chosen so that thermodynamical 
\eqref{td-300} and \eqref{td-30} are satisfied. Using \eqref{res-0} and \eqref{res-10} we obtain
\begin{equation} \label{res-3}
u(x,t)=\frac{1}{\pi} \int_{0}^{t} \Big[ \int_{0}^{\infty}
\exp \big[ \tau\imp M(i\tau) x \big] \cos \big[ \tau (\theta -\rep M(i\tau) x) \big] \,d\tau \Big] \,d\theta.
\end{equation}
For $x=0$ from \eqref{res-3} we obtain
$$
u(0,t)=\frac{1}{\pi} \int_{0}^{t} \int_{0}^{\infty} \cos (\tau \theta) \,d\tau \,d\theta.
$$
By using the fact that Dirac $\delta$ function may be approximated as
$$
\delta (\xi) = \frac{1}{\pi} \lim_{\nu \to 0} \int_{0}^{1/\nu} \cos (\tau \theta) \,d\tau,
$$
we obtain
$$
u(0,t)=1.
$$
Thus the boundary condition $u(0,t)=H(t)$ is satisfied.

\section{Conclusion}

In this work we proposed a constitutive equation for viscoelastic body of generalized 
Zener type that includes fractional derivatives of stress and strain of real \textit{and complex}
order. With such constitutive equation, the initial-boundary value problem that generalizes 
the classical wave equation is given by \eqref{cfZwe}-\eqref{bc1}. Note
that for the case $a_{1}=b_{1}=a_{2}=b_{2}=0$ the constitutive equation \eqref{1a}$_{2}$ 
becomes Hooke's law, and \eqref{cfZwe}-\eqref{bc1} reduces to an initial-boundary value
problem for the classical wave equation.

The results of this paper may be summarized as follows:
\begin{enumerate}
\item We formulated initial-boundary value problem for the generalized wave equation in the 
viscoelastic body described by fractional derivatives of real and complex order in the form
(\ref{cfZwe})-(\ref{l-1}).
\item We determined restrictions on the coefficients from the dissipativity conditions in the 
form \eqref{td-1}, \eqref{td-300} and \eqref{td-30}, and later in a strong form, for the sake of solvability. 
We concluded that dissipativity conditions, 
that are consequences of the Second Law of Thermodynamics for isothermal
deformation in a strong form \eqref{td-1}, \eqref{td-300}$_s$ and \eqref{td-30}$_s$, 
guarantee the solvability of the constitutive equation \eqref{1a}$_{2}$ for $\sigma$ and $\eps$.
\item We presented the solution to \eqref{cfZwe}-\eqref{bc1} in the form of \eqref{res-0}.
\item We analyzed two specific examples. In the first example the solution is given by 
\eqref{res-11}. For the Dirac delta impulse as the boundary condition, the solution shows 
a pulse behavior that dissipates with time. In calculating the integral in \eqref{res-11} 
there was observed oscillation type behavior for small times. We attribute this behavior 
to the imaginary order term in the derivation since for small $p$ the influence of oscilating factor in the integrals $I_2, I_6$
is evident. 
\item Further study is needed to examine the influence of parameters of the model and 
order of derivatives on the properties of the solution.
\end{enumerate}

 \section*{Acknowledgements}

 This work is supported by
 Projects 174005 and 174024 of the Serbian Ministry of Science.



\begin{thebibliography}{10}

\bibitem{AmendolaFabrizioGolden}
{Amendola, G., Fabrizio, M., Golden, J.~M.}
\newblock {\em Thermodynamics of Materials with Memory}.
\newblock Springer, New York, 2012.

\bibitem{Atanackovic02}
{Atanackovi\' c, T.M.}
\newblock A modified {Z}ener model of viscoelastic body.
\newblock {\em Contin. Mech. Thermodyn.}, {\bf 14}:137--148, 2002.

\bibitem{AJKPZ15}
{Atanackovi\' c, T.M., Janev, M., Konjik, S., Pilipovi\' c, S., Zorica, D.}
\newblock Vibrations of an elastic rod on a viscoelastic foundation of complex
  fractional kelvin-voigt type.
\newblock {\em Meccanica}, {\bf 50}(7):685--704, 2015.

\bibitem{AKOZ-thermodynamical}
{Atanackovi\' c, T.M., Konjik, S., Oparnica, Lj., Zorica, D.}
\newblock Thermodynamical restrictions and wave propagation for a class of
  fractional order viscoelastic rods.
\newblock {\em Abstr. Appl. Anal.}, {\bf 2011}:975694(32pp), 2011.

\bibitem{AKPZ-complexce}
{Atanackovi\' c, T.M., Konjik, S., Pilipovi\' c, S., Zorica, D.}
\newblock Complex order fractional derivatives in viscoelasticity.
\newblock {\em Mech. Time-Depend. Mater.}, 2016.

\bibitem{APZ10-creep}
{Atanackovi\' c, T.M., Pilipovi\' c, S., Zorica, D.}
\newblock Distributed-order fractional wave equation on a finite domain.
  {C}reep and forced oscillations of a rod.
\newblock {\em Contin.\ Mech.\ Thermodyn.}, {\bf 23}:305--318, 2011.

\bibitem{APZ10-stress}
{Atanackovi\' c, T.M., Pilipovi\' c, S., Zorica, D.}
\newblock Distributed-order fractional wave equation on a finite domain.
  {S}tress relaxation in a rod.
\newblock {\em Int.\ J.\ Eng.\ Sci.}, {\bf 49}:175--190, 2011.

\bibitem{APSZ-book-vol1}
{Atanackovi\'c, T.\ M., Pilipovi\'c, S., Stankovi\'c, B., Zorica, D.}
\newblock {\em Fractional Calculus with Applications in Mechanics: Vibrations
  and Diffusion Processes}.
\newblock Wiley-ISTE, London, 2014.

\bibitem{APSZ-book-vol2}
{Atanackovi\'c, T.\ M., Pilipovi\'c, S., Stankovi\'c, B., Zorica, D.}
\newblock {\em Fractional Calculus with Applications in Mechanics: Wave
  Propagation, Impact and Variational Principles}.
\newblock Wiley-ISTE, London, 2014.

\bibitem{BagleyTorvik}
{Bagley, R. L., Torvik, P. J.}
\newblock On the fractional calculus model of viscoelastic behavior.
\newblock {\em J.\ Rheology}, {\bf 30}(1):133--155, 1986.

\bibitem{Doetsch-Handbuch-1}
{Doetsch, G.}
\newblock {\em Handbuch der Laplace-Transformationen I}.
\newblock Birkh\"{a}user, Basel, 1950.

\bibitem{Franchi}
{Franchi, F., Lazzari, B., Nibbi, R.}
\newblock Mathematical models for the non-isothermal {J}ohnson-{S}egalman
  viscoelasticity in porous media: stability and wave propagation.
\newblock {\em Math. Meth. Appl. Sci.}, {\bf 38}:4075--4087, 2015.

\bibitem{Hanyga-07}
{Hanyga, A.}
\newblock Fractional-order relaxation laws in non-linear viscoelasticity.
\newblock {\em Contin. Mech. Thermodyn.}, {\bf 19}:25--36, 2007.

\bibitem{Hanyga-16}
{Hanyga, A.}
\newblock Wave propagation in anisotropic viscoelasticity.
\newblock {\em J Elasticity}, {\bf 122}(2):231--254, 2016.

\bibitem{KOZ10-WZ}
{Konjik, S., Oparnica, Lj., Zorica, D.}
\newblock Waves in fractional {Z}ener type viscoelastic media.
\newblock {\em J.\ Math.\ Anal.\ Appl.}, {\bf 365}(1):259--268, 2010.

\bibitem{Lion}
{Lion, A.}
\newblock On the thermodynamics of fractional damping elements.
\newblock {\em Contin. Mech. Thermodyn.}, {\bf 9}:83--96, 1997.

\bibitem{Love-complex}
{Love, E. R.}
\newblock Fractional derivatives of imaginary order.
\newblock {\em J. London Math. Soc.}, {\bf 2-3}(2):241--259, 1971.

\bibitem{Mainardi10}
{Mainardi, F.}
\newblock {\em Fractional Calculus and Waves in Linear Viscoelasticity}.
\newblock Imperial College Press, London, 2010.

\bibitem{MakrisC-91}
{Makris, N., Constantinou, M.}
\newblock Fractional-derivative {M}axwell model for viscous dampers.
\newblock {\em J Struct. Eng.}, {\bf 117}:2708--2724, 1991.

\bibitem{MakrisC-92}
{Makris, N., Constantinou, M.}
\newblock Spring-viscous damper systems for combined seismic and vibration
  isolation.
\newblock {\em Earthq. Eng. Struct. Dynam.}, {\bf 21}:649--664, 1992.

\bibitem{MakrisC-93}
{Makris, N., Constantinou, M.}
\newblock Models of viscoelasticity with complex-order derivatives.
\newblock {\em J. Eng. Mech.}, {\bf 119}(7):1453--1464, 1993.

\bibitem{Podlubny}
{Podlubny, I.}
\newblock {\em Fractional Differential Equations}, volume 198 of {\em
  Mathematics in Science and Engineering}.
\newblock Academic Press, San Diego, 1999.

\bibitem{SKM}
{Samko, S. G., Kilbas, A. A., Marichev, O. I.}
\newblock {\em Fractional Integrals and Derivatives - Theory and Applications}.
\newblock Gordon and Breach Science Publishers, Amsterdam, 1993.

\bibitem{Tamaogi}
{Tamaogi, T., Sogabe, Y.}
\newblock Longitudinal wave propagation including high frequency component in
  viscoelastic bars.
\newblock In {Song, B., Lamberson, L., Casem, D., Kimberley, J.}, editor, {\em
  Dynamic Behavior of Materials, Volume 1 Proceedings of the 2015 Annual
  Conference on Experimental and Applied Mechanics}, pages 75--80. Springer,
  2016.

\bibitem{Wang-16}
{Wang, Y.}
\newblock Generalized viscoelastic wave equation.
\newblock {\em Geophys. J. Int.}, {\bf 204}:1216--1221, 2016.

\end{thebibliography}
 \end{document}